\documentclass[12pt]{article}
\font\sdopp=msbm10
\def\CI {\sdopp {\hbox{C}}}
\date{07/12/2003}
\title{
Alguns fenomens 
de continuaci\'o
anal\'\i tica\\
en una variable complexa
}
\author{Claudi Meneghin}

\newtheorem{definition}{Definition}

\font\sdopp=msbm10

\def\ERRE {\sdopp {\hbox{R}}}
\def\CI {\sdopp {\hbox{C}}}
\def\TI {\sdopp {\hbox{T}}}
\def\DI {\sdopp {\hbox{D}}}
\def\ENNE {\sdopp {\hbox{N}}}
\def\ZETA{\sdopp {\hbox{Z}}}

\def\ce {\sdopp {\hbox{C}}}

\def\ZETA{\sdopp {\hbox{Z}}}

\def\ii{i} 
\def\e{\hbox{\boldmath{}$e$\unboldmath}}

\def\id{\hbox{\boldmath{}$id$\unboldmath}}

\newtheorem{theoreme}[definition]
{Teorema}

\begin{document}
\maketitle
\bibliographystyle{plain} 
\parindent=8pt
%
\def\Ti {\sdopp {\hbox{T}}}
\def\IM{\hbox{\boldmath{}$i$\unboldmath}} 

\def\Ch{\hbox{\rm Ch}}
\def\CIRC{\mathop{\tt o}\limits}
\def\BBB{\sl}


\font\sdopp=msbm10
\def\de {\sdopp {\hbox{D}}}

\def\M{\hbox{\tt\large M}}
\def\N{\hbox{\tt\large N}}
\def\T{\hbox{\tt\large T}}

\def\e{\hbox{\boldmath{}$e$\unboldmath}} 

\def\labelle #1{\label{#1}}
\begin{abstract}
In this expository paper we examinate some
phenomena arising when
a holomorphic germ is
analytically continued.
\end{abstract}

\section{Introducci\'o}
Es consideri un punt $p$
del pla complex i
la 
s\`erie de
pot\`encies en $z-p$:
$$
\alpha_0+\alpha_1(z-p)+
\alpha_2(z-p)^2+
\alpha_3(z-p)^3+...
$$

Aquesta s\`erie convergeix en un cert cercle $C_1$ de centre $p$ i doncs hi definiex 
una funci\'o
holomorfa $f$; escrivem $ f_p  $
per a posar en evid\`encia  el punt de 
desenvolupament.

Considerem un punt $q\in C_1$
i desenvolupem $f$ en 
s\`erie de pot\`encies 
de $z-q$:
$$
f_q(z)=\beta_0+\beta_1(z-q)+
\beta_2(z-q)^2+
\beta_3(z-q)^3+...
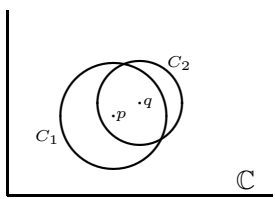
$$

\begin{tabular}{cc}
\thicklines
\begin{picture}(100,40)(20,30)
\put(10,0){\line(1,0){100}}
\put(10,0){\line(0,1){70}}
\put(100,6){\makebox(0,0){\tiny$\CI$}}
\put(50,30){\circle{40}}
\put(60,35){\circle{30}}
\put(50,30){\circle*{1}}
\put(60,35){\circle*{1}}
\put(53,30){\makebox(0,0){\tiny$p$}}
\put(63,35){\makebox(0,0){\tiny$q$}}
\put(25,22){\makebox(0,0){\tiny$C_1$}}
\put(75,50){\makebox(0,0){\tiny$C_{2}$}}
\end{picture}
&
\begin{minipage}{262pt}
Si \'es cas que 
 el cercle de 
converg\`encia $C_2$ d'aquesta darrera
s\`erie 
no sigui continugut en $C_1$,
  hom
ha de fet  obtingut una coneixen\c ca
m\'es ampla de $f$, mitjan\c cant la definici\'o:
$$
f(z):=\cases{
f_p(z)
& si $z\in C_1$\cr
f_q(z)
& si $z\in C_2$
}.
$$
\end{minipage}
\end{tabular}
\vskip0,3truecm

Aquesta definici\'o  \'es 
b\'e posada, perqu\`e $ z\in C_1\cap C_2 \Rightarrow  f_p(z)
=
f_q(z)
$.

Direm que l'extenci\'o 
de $ f  $ a 
$C_1\cup C_2   $ 
aix\'\i\   obtinguda \'es una 
{\it continuaci\'o anal\'\i tica }
(o tamb\'e un {\it prolungament anal\'\i tic})
de $ f_p :C_1\rightarrow \CI $;
direm tamb\'e que $ f_q :C_2\rightarrow \CI $
\'es una 
continuaci\'o anal\'\i tica de
$ f_p :C_1\rightarrow \CI $ i viceversa.
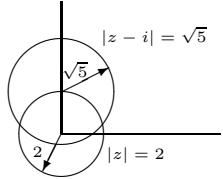
\begin{tabular}{cc}
\begin{picture}(86,50)(20,30)
\put(50,20){\line(1,0){60}}
\put(50,20){\line(0,1){50}}
\put(50,20){\circle{32}}
\put(50,36){\circle{40}}
\put(78,12){\makebox(0,0){\tiny$\vert z\vert={2}$}}
\put(85,57){\makebox(0,0){\tiny$\vert z- i\vert=\sqrt{5}$}}
\put(50,20){\vector(-1,-2){7}}
\put(50,36){\vector(2,1){18}}
\put(41,13){\makebox(0,0){\tiny $2$}}
\put(55,44){\makebox(0,0){\tiny 
$\sqrt{5}$}}
\end{picture}
&
\begin{minipage}{286pt}
Per exemple, es pot senzillament
veure que les dues 
s\`eries de pot\`encies 
$\displaystyle\sum_{n=0}^{\infty}\frac{z^n}
{2^{n+1}}\ (\vert z     \vert<2)$ i 
$\displaystyle\sum_{n=0}^{\infty}\frac{(z-i)
^n}
{(2-i)^{n+1}}\ \vert z-i     \vert<\sqrt{5}$ s\'on 
cadascuna continuaci\'o anal\'\i tica
de l'altra.
\end{minipage}
\end{tabular}

\vskip0,5truecm
 Notem que totes dues representen 
la funci\'o $z\mapsto
 1/(2-z)$. M\'es in general,
si \'es cas que 
$f$, definida
a priori dins un conjunt 
obert
$U\subset\CI$, es pugui
restringir a un conjunt obert
$V\subset U$ i successivament
$f\vert_V$ pugui \'esser
prolongada
a un conjunt obert 
$W\not\subset U$, direm que 
la nova
funci\'o obtenida es una
'continuaci\'o anal\'\i tica'
de $f$. 

Aquest article t\'e a veure 
amb alguns 
fenomens interessants 
que appereixen de manera
natural
en un tal context:
il$\cdot$lustrarem 
la formaci\'o de fronteres 
naturals per al prolungaments
anal\'\i tics,
el fenomen
de
la resurg\`encia i un
fenomen de idiosincrasia
amb la noci\'o de revestiment
topol\`ogic.

L'escolliment dels arguments
 \'es 
una conseq\"u\`encia del gust de
l'autor.

\subsection{ Les definicions b\`asiques}
Un {\sf  element de funci\'o
holomorfa}
\'es un parell
$\left(U,f    \right)$, 
on
$U$ \'es
un conjunt obert
a conexi\'o simple
del
pla complex,
$f:U\rightarrow \CI$ una funci\'o
 holomorfa  definida en $U$, que pren valors
en $\ce$.
Dos elements
$\left(U,f    \right)$ i $\left(V,g    \right)$ 
s\'on {\sf conectables} si
existeix
 una successi\'o finita 
$ \left\{(U_j,f_j)
    \right\}_{j=0,....,n}$,
tal que
$\left(U_0,f_0    \right)=\left(U,f    \right)$, $\left(U_n,f_n    \right)=\left(V,g    \right)$
i,
per a tot 
$j=0,....,n-1$,
$$
\cases
{
U_j\cap U_{j+1}\not= \emptyset,\cr 
f_{j+1}\vert_{U_j\cap U_{j+1}}=f_{j}
\vert_{U_j\cap U_{j+1}}.
}
$$
Direm que
$\  \{(U_i,f_i)\}_{i=0...n}\   $ 
\'es una {\bf continuaci\'o anal\'\i tica}
de $ (U,f)  $  (o de $ (V,g)  $).
Direm tamb\'e, si no hi ha 
possibilitat de confusi\'o,
que cada element

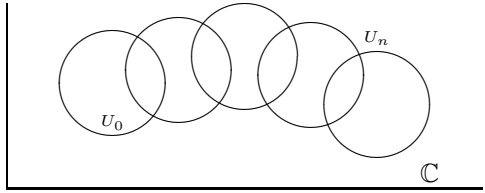
\begin{tabular}{ccc}
\begin{picture}(100,50)(20,30)
\put(10,0){\line(1,0){180}}
\put(10,0){\line(0,1){70}}
\put(170,6){\makebox(0,0){\tiny$\CI$}}
\put(50,25){\makebox(0,0){\tiny$U_0$}}
\put(150,57){\makebox(0,0){\tiny$U_n$}}
\put(50,40){\circle{40}}
\put(75,45){\circle{40}}
\put(100,50){\circle{40}}
\put(125,43){\circle{40}}
\put(150,32){\circle{40}}
\end{picture}
&
\begin{minipage}{90pt}
$$
\ 
$$
\end{minipage}
&
\begin{minipage}{162pt}
\'es una continuaci\'o
anal\'\i tica de $(U,f)$ (o de $ (V,g)  $).
Els elements $ \left\{(U_j,f_j)
    \right\}_{j=0,....,n}$
es diran 
{\sf enlla\c cats}.

Una continuaci\'o
anal\'\i tica $\  \{(U_i,f_i)\}_{i=0...n}\   $ 
es dir\`a una {\sf continuaci\'o
anal\'\i tica al llarg d'un}

\end{minipage}
\end{tabular}
{\sf cam\'\i\ }  
$\gamma:[0,1]\rightarrow\ce$
(que, per a 
senzillesa, suposem $C^1$ a trets)
si esisteix una partici\'o $\{I_i\}_{i=0...n}$
de $[0,1]$ tal que $\gamma (I_i)\subset U_i$ i
$\gamma\vert_{(I_i)}$ \'es una aplicaci\'o
injectiva per a cada $i=0...n$.

Cal sens dubte recordar que 
la continuaci\'o
anal\'\i tica al llarg d'un 
cam\'\i\  tancat 
no conserva pas, en general,  
els valors de la funci\'o
en un entorn
del punt de partida:
 es tingui en compte,
per exemple, 
la determinaci\'o $\varphi$
de
la funci\'o 'arrel 
quadrada complexa',
en un entorn  de $1$,
tal que $\varphi(1)=1$.

Es pot veure $\varphi$, 
en coordenades polars,
com a l'aplicaci\'o que
envia 
$\varrho \exp(i\vartheta)$ 
vers
$\sqrt{\varrho}\exp(i\vartheta/2)$, on
$\sqrt{\ }$ indica l'operaci\'o d'arrel quadrada
real positiva. 
Intu\"\i tivament,
continuem
$\varphi$
al llarg de la 
circumfer\`encia unitat:
despr\'es
una volta compleda, \'es a dir
un increment de $\vartheta   $
igual  a $ 2\pi  $, 
obtenim 
un nou
element de funci\'o holomorfa 
$\psi$ en un entorn de $1$,
que ha 
redu\"\i t a meitat 
l'increment de de $2\pi$ 
l'argument de $\zeta$.

Doncs, $\arg(\psi(\zeta))
=\arg(\varphi(\zeta))+\pi$,
\'es a dir $\varphi=-\psi$.
Naturalment, una altra volta
de $2\pi$ ens porta 
de bell nou a l'element 
de partida $\varphi$.

\subsection{Alguns enunciats
cl\`assics\labelle{intro}}

Per a acabar aquesta 
introducci\'o,  recordem
alguns enunciats
estandard d'an\`alisi complexa i topologia
general, que 
seran utilitzats 
m\'es enll\`a.

\begin{theoreme}\labelle{successio}
{\tt (Teorema d'unicitat),\cite{ahlfors} 
} Siguin $f $ i $ g  $ dues funcions meromorfes
en un conjunt obert $ U\subset\CI  $: suposem $ f=g  $ en un conjunt amb un punt
d'acumulaci\'o $ p\in U  $: llavors $ f\equiv
g  $ en $ U  $.
\end{theoreme}
\begin{theoreme}\labelle{mapa}
{\tt (Teorema de l'aplicaci\'o de Riemann)
\cite{ahlfors} }
Sigui $ U\subset\CI  $, $ U\not=\CI  $,
un conjunt obert, a connexi\'o simple:
llavors $ U  $ \'es biholom\`orfic a $ \DI  $.
\end{theoreme}
\begin{theoreme}\labelle{koebe}
{\tt (Teorema-
$ 1/4  $ de Koebe) \cite{pommerenke}}
Si $ f:\DI\rightarrow \CI  $
 \'es una mapa conforme, llavors
$$
\displaystyle
\frac{1}{4}
\left(1-\vert z     \vert^2    \right)
\vert f^{\prime}(z)     \vert
\leq dist\left(f(z),\partial f(\DI)    \right)
\leq
\left(1-\vert z     \vert^2    \right)
\vert f^{\prime}(z)     \vert
$$
per a tot $ z\in\DI  $.
\end{theoreme}

Siguin $ X  $ i $ Y  $ uns espais topol\`ogics: 
una aplicaci\'o continua surjectiva
$p:Y\rightarrow X$
 \'es un 
{\bf  revestiment topol\`ogic}
si cada punt 
$x\in X$
t\'e un entorn obert $ {\cal U}  $
tal que la restricci\'o de
$p$ a cada component $ {\cal V}_i  $
de
$p^{-1}({\cal U}  )$ \'es un 
homeomorfisme  de $ {\cal V}_i  $ sobre
${\cal U}$. 

Tamb\'e recordem que 
una aplicaci\'o continua 
 $p:Y\rightarrow X$
t\'e la 
{\bf proprietat de l'elevament de les corbes
} si, per a cada corba
$\gamma:I\rightarrow X$ i cada 
$y\in p^{-1}(\gamma(0))   $
existeix una corba 
$\widetilde \gamma:I\rightarrow Y   $
tal que $p\circ\widetilde\gamma=
\gamma$ i $\widetilde \gamma(0)=y$.

El resultat seg\"uent \'es estandard:
(vegeu per exemple \cite{klaus}, secci\'o  9.3
i \cite{forster}, teorema 4.19) 
\begin{theoreme}
Un homeomorfisme local
surjectiu 
entre dos
espais topol\`ogics
\'es un revestiment topol\`ogic
si i solament si t\'e la proprietat 
de l'elevament de les corbes.
\end{theoreme}

\section{Formaci\'o de fronteres naturals}
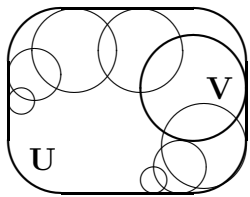
\begin{tabular}{cc}
\begin{picture}(110,45)(0,35)
{
\thicklines
\put(50,40){\oval(90,70)}
\put(75,45){\circle{40}}
}
\put(18,18){\makebox(0,0){\bf U}}
\put(79,23){\circle{30}}
\put(55,59){\circle{30}}
\put(30,59){\circle{30}}
\put(15,50){\circle{20}}
\put(10,40){\circle{10}}
\put(70,15){\circle{20}}
\put(60,10){\circle{10}}
\put(85,45){\makebox(0,0){\bf V}}
\end{picture}
&
\begin{minipage}{262pt}
Es consideri un element 
de funci\'o holomorfa $(U,f)   $:
pot succeir que,
per a cada restricci\'o $ (V,g )  $ 
de $(U,f)     $  (\'es a dir,
$ V\subset U  $ i $ g=f\vert_V  $)
no existeixi cap continuaci\'o anal\'\i tica
$ (W,h)  $ de  $ (V,g )  $  tal que 
$
W\cap U\not\subset U$.
Si \'es cas, direm que $ \partial U$
\'es
una {\bf frontera natural} per a 
l'element $(U,f)   $.
\end{minipage}
\end{tabular}
%
%
%
Considerem per exemple
la s\'erie de pot\`encies 
$\sum_{n=0}^{\infty}z^{2^n}=1+z^2+z^4+
z^8+...$:
gracies al teorema de Cauchy-Hadamard
ella convergeix dins el disc $\vert z\vert<1$,
i doncs hi definieix una funci\'o 
holomorfa $ h  $.
De m\'es, 
$ h(z)\to\infty  $
llavors que $z\to 1$ al llarg de l'eix real.
Puix que
$
h(z^2)=1+z^4+z^8+
z^{16}+...=h(z)-z^2
$ hom ha 
$
\lim_{z\to{-1},z\in\ERRE}
h(z)
=$
$
\lim_{z\to{-1},z\in\ERRE}
\left(z^2+h(z^2)    \right)
=\infty
$.

De la mateixa manera, $h(z)=z^2+z^4+h(z^4)$,
doncs $h\to\infty$ 
llavors que $z\to\pm i$ al llarg de l'eix 
imaginari; de manera general,
$
h(z)=z^2+...+z^{2^n}+h(z^{2^n})
$,
 per a tot nombre natural $n   $, doncs
$h\to\infty$ 
llavors que $z\to \exp({2k\pi\ii/2^n})$ 
al llarg d'un radi del disc.

El conjunt dels punts de la forma 
$\exp({2k\pi\ii/2^n}),\ k,n\in\ZETA$ 
\'es dens dins el cercle
$\TI=\{\vert z\vert=1\}   $, doncs 
$ h  $ no admet cap continuaci\'o anal\'\i tica
a algun punt d'aquesta corba: ella \'es doncs 
una frontera natural. 

Observem que $h$ pot tampoc 
ser continuada als punts de $ \TI  $ com a
funci\'o meromorfa, perqu\`e, en aquest cas,
$1/h   $ 
s'anul$\cdot$laria en un conjunt amb un punt d'acumulaci\'o  i seria doncs
identicament zero (vegeu th.\ref{successio})  
(Pel que concerneix aquest exemple,
el lector podr\`a tamb\'e consultar,
per exemple,
\cite{chabat} p.129).

\section{La resurg\`encia}
En A.Hurwitz va plantejar, en el seu quadern,
a la data del 6 desembre 1918, 
la demanda si fou possible
 que una 
s\`erie de pot\`encies 
\begin{equation}
h(\xi)=\sum_{k=0}^{\infty}a_k(\xi-\xi_0)^k,
\labelle{hur_1}
\end{equation}
representant una funci\'o diferent de
$ \xi\mapsto c\e^{\xi}  $,
admet\'es continuaci\'o anal\'\i tica al llarg 
d'un cam\'\i \ tancat $\gamma$ al voltant de
$ \xi_0  $ i, a la fi
de la continuaci\'o, prengu\'es la forma
\begin{equation}
\sum_{k=1}^{\infty}k a_k(\xi-\xi_0)
^{k-1}=h^{\prime}(\xi).
\labelle{hur_2}
\end{equation}

\subsection{La soluci\'o de Lewy}
En H.Lewy
va donar una soluci\'o del problema
(\ref{hur_1})/(\ref{hur_2})
,
que presentem aqu\'\i\  
en una forma lleugerament modificada.

Es consideri la funci\'o:
$$
h(z)=\int_{\ERRE^+}
 \exp
\left[
-zt-(\log t)^2/4\pi i 
\right]\, dt;
$$
$ h  $ \'es holomorfa per 
$\Re(z)>0 $
i pot ser continuada 
anal\'\i ticament als semiplans 
$\Re(z \e^{- i\vartheta})>0\ (\vartheta
\in\ERRE^+)$,
de la manera seg\"uent:
sigui $ N\in\ENNE  $ tal que
$0<\vartheta/N<\pi/2  $ i fem
$\eta:= \vartheta/N  $.

Escrivem, per a $ z\in
\Re(z \e^{- i\eta})>0\bigcup
\Re(z)>0 $,
\begin{eqnarray}
h(z) &=&
 \int_{\ERRE^+}
 \exp
\left[
z
\e^{- i\eta}
\e^{i\eta}
t-
\frac{\log(
\e^{- i\eta}
\e^{i\eta}
t)^2}{4\pi i  }
\right]
\, dt\nonumber\\
&=&
\int_{\e^{i\eta}\ERRE^+}
 \exp
\left[-z\e^{- i\eta}
u-
\displaystyle
\frac{(\log(u)-i\eta)^2}{4\pi i }
\right]
\e^{- i\eta}
\, du\nonumber\\
&=&
\lim_{R \to\infty}
\left\{
\int_0^{R}
 \exp
\left[-z\e^{- i\eta}
u-
\displaystyle
\frac{(\log(u)-i\eta)^2}{4\pi i }
\right]
\e^{- i\eta}
\, du+\right.\nonumber\\
&\ &\qquad \left. + 
\int_{\gamma_R}
 \exp
\left[-z\e^{- i\eta}
u-
\displaystyle
\frac{(\log(u)-i\eta)^2}{4\pi i }
\right]
\e^{- i\eta}
\, du\right\} 
\labelle{integralcurv}.
\end{eqnarray}

La integral en (\ref{integralcurv}),
que anominem $ I_2  $,
 ha de ser 
calculada sobre la corba 
$\gamma_R:
[0,1]\rightarrow\CI$ definida posant
$ \gamma (t):= R\e^{it\eta} $.

Hom ha $ I_2\leq C_1 R^{\alpha}\e^{-C_2R}  $
per a unes constantes reals positives 
$C_1$, $ {C_2}  $ i ${\alpha}     $, doncs 
$ I_2  $
tendeix a
$ 0  $ quan $ R\to\infty$.

Aix\'\i\ 
per a
 $ z\in
\{\Re(z \e^{- i\eta})>0     \}\bigcap
\{\Re(z)>0     \} $ hom ha 
\begin{eqnarray*}
h(z)
&=&
\int_{\ERRE^+}
 \exp
\left[-z\e^{- i\eta}
u-
\displaystyle
\frac{(\log(u)-i\eta)^2}{4\pi i }
\right]
\e^{- i\eta}
\, du;
\end{eqnarray*}
per\`o
aquesta darrera integral convergeix en
$\Re(z \e^{- i\eta})>0$
i doncs hi defineix
una continuaci\'o  anal\'\i tica de $ h  $.
Repetem el procediment $ N  $
vegades: aix\`o ens dona finalment una continuaci\'o 
anal\'\i tica de $ h  $ al semipl\`a
$\Re(z \e^{- i\vartheta})>0$; doncs $ h  $
pot ser continuada anal\'\i ticament
 a tot punt
 $p\in\CI\setminus\{0     \}  $.

Finalment, si fem la continuaci\'o  anal\'\i tica al llarg del cam\'\i \ $ \vert z\vert=1,
0\leq\arg(z)\leq 2\pi   $, obtenim, designant
$\widehat h   $ l'element de funci\'o holomorfa
obtingut (en un entorn de $ z=1  $) despr\'es una volta completa,
\begin{eqnarray*}
\widehat h(z) &=&
 \int_{\ERRE^+}
 \exp
\left[-\e^{2\pi i}z
t-(\log t+2\pi i)^2/4\pi i    \right]\, dt\\
&=&
\int_{\ERRE^+}
 \exp
\left[-zt-
\displaystyle
\frac{(\log t)^2-4\pi ^2+4\pi i\log t}{4\pi i }
 \right]\, dt\\
&=&
\int_{\ERRE^+}
 \exp
\left[
\displaystyle
-zt
-\e^{2\pi i}
t-(\log t)^2/4\pi i
-
\pi i+ \log t
 \right]\, dt\\
&=&
\int_{\ERRE^+}
(-t) 
\exp
\left[-zt-(\log t)^2/4\pi i    \right]\, dt\\
&=&
\Big.
h^{\prime}(z).
\end{eqnarray*}

\subsection{El m\`etode de la integral de Laplace}

Es consideri altra vegada el problema 
de Hurwitz 
(\ref{hur_1})/(\ref{hur_2})
i es suposi $\xi_0\not=0$; el cam\'\i \
$\gamma$ sigui el cercle $\vert\xi\vert=\vert\xi_0\vert$
amb l'orientaci\'o positiva.

Fem $\xi=\exp(2\pi i z)$ i 
$h(\exp(2\pi i z))=G(z)$: aix\`o transforma el problema 
en el seg\"uent:
existeixen solucions de l'equaci\'o 
\begin{equation}
G(z+1)=(1/2\pi i)\exp(-2\pi iz)G^{\prime}(z)
\labelle{quater}
\end{equation}
que siguin 
anal\'\i tiques in una tira horitzontal
del pla complex i differentes de 
$c \exp(\exp(2\pi i z))$?

Aquest darrer problema va \'esser 
analitzat, de manera m\'es general, per 
A.Naftalevich (vegeu \cite{naftalevich})
i Berenstein i Sebbar (vegeu \cite{berenseb})
.

Presentem aqu\'\i \ el m\`etode de la integral
de Laplace, descrit en \cite{berenseb}, sec. 2,
que produeix solucions enteres
de (\ref{quater}): aqueste m\`etode
consisteix a 
 trobar
una condici\'o necessaria
per a la soluci\'o de (\ref{quater})
i a demostrar que ella \'es tamb\'e
suficient.

Suposem que l'equaci\'o (\ref{quater}) tingui
solucions enteres de la  forma 
$ G(z)=\int_C \e^{-2i\pi uz} \varphi(u)\,du $,
per una funci\'o  holomorfa 
$ \varphi  $ 
i  un cam\'\i \ d'integraci\'o
$ C  $ oportuns.

Hom ha
\begin{eqnarray}
G(z+1)&=&\int_C \e^{-2i\pi u(z+1) }\varphi(u)\,du
\nonumber\\
&=&\int_{C-1} \e^{-2i\pi (u+1)(z+1) }\varphi(u+1)\,du
\nonumber\\
&=&\int_{C-1} \e^{-2i\pi (u+1)z}
\e^{-2i\pi u}
\varphi(u+1)\,du
\labelle{camia}\\
&=&\int_{C} \e^{-2i\pi (u+1)z}
\e^{-2i\pi u}
\varphi(u+1)\,du.
\labelle{camib}
\end{eqnarray}

La deducci\'o de (\ref{camib}) a partir de 
(\ref{camia}) es justificada per l'hip\`otesi
que $ \varphi  $ sigui entera i $ C  $ oport\'u.

Hom ha tamb\'e
$$
\frac{1}{2\pi i}
\e^{-2\pi iz}G^{\prime}(z)
=
-\int_C u\e^{-2i\pi (u+1)z} \varphi(u)\,du 
=
-\int_{C} (u-1)\e^{-2i\pi uz} \varphi(u-1)\,du. 
$$

Aix\'\i , el problema ha estat reduit a trovar
solucions holomorfes (en una tira horitzontal) de l'equaci\'o 
$ \e^{-2i\pi u} \varphi(u+1)=-u\varphi(u)  $:
si fem $\varphi(u)=\e^{i\pi u^2}\Phi(u)$ aix\`o \'es el mateix que la cerca d'una funci\'o  holomorfa $ \Phi  $ tal que 
$\Phi(u+1)=u \Phi(u)    $.

Naturalment, la funci\'o $\Gamma   $ 
de Euler \'es una soluci\'o, doncs hem arribat
a la soluci\'o  formal de (\ref{quater}):
\begin{equation}
\labelle{formal}
G(z):=\int_{C}\e^{-2i\pi uz+i\pi u^2}\Gamma  (u)\, du.
\end{equation}
 
Sigui $ C  $ la recta real $ t\mapsto 1+(1+i)t  $:
la representaci\'o asimpt\`otica de Stirling:
$$
\Gamma  (u)=\exp\left[
\left(u-1/2    \right)\hbox{\tt LOG}(u)
-u+\frac{1}{2}\log(2\pi )
    \right]\left(1+O(1/\vert u\vert)    \right),
\ \vert \arg(u)\vert<\pi, 
$$
(on \hbox{\tt LOG} \'es la determinaci\'o 
principal del logaritme)
implica que la integral en (\ref{formal})
convergeix, perqu\`e
\begin{equation}
\e^{-2i\pi u z+i\pi u^2}\Gamma  (u)
\vert_{u=t(1+i)+1}
=
O\left[\e^{-2\pi t^2 + A\vert t
     \vert}\log(B\vert t     \vert)    \right]
\labelle{grando}
\end{equation}
per a dues constantes $ A,B\in\ERRE^+  $
oportunes quan $ \vert t     \vert\to\infty$;
aix\'\i\  $ G  $ \'es una funci\'o 
entera.

Es pot tamb\'e veure que el cam\'\i\ 
$ C  $ escollit \'es oport\'u: puix que 
\begin{eqnarray*}
&\ &\int_{C-1} \e^{-2i\pi (u+1)z}
\e^{-2i\pi u}
\varphi(u+1)\,du
\\
&=&\lim_{t\to\infty}
\int_{D_{-t}}+\int_{C_{t}}+\int_{D_{t}}\\
&=& I_{1t}+I_{2t}+I_{3t}, 
\end{eqnarray*}
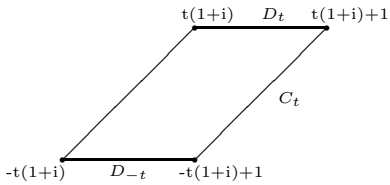
\begin{tabular}{cc}
\begin{picture}(130,40)(0,30)
\put(10,10){\line(1,0){50}}
\put(10,10){\circle*{2}}
\put(10,10){\line(1,1){50}}
\put(60,10){\line(1,1){50}}
\put(60,10){\circle*{2}}
\put(60,60){\line(1,0){50}}
\put(60,60){\circle*{2}}
\put(110,60){\circle*{2}}
\put(0,5){\makebox(0,0){\hbox{\tiny -t(1+i)}}}
\put(35,5){\makebox(0,0){\hbox{\tiny $D_{-t}$}}}
\put(70,5){\makebox(0,0){\hbox{\tiny -t(1+i)+1}}}
\put(65,65){\makebox(0,0){\hbox{\tiny t(1+i)}}}
\put(90,65){\makebox(0,0){\hbox{\tiny $D_{t}$}}}
\put(120,65){\makebox(0,0){\hbox{\tiny t(1+i)+1 }}}
\put(96,33){\makebox(0,0){\hbox{\tiny $C_{t}$}}}
\end{picture}
&
\begin{minipage}{242pt}
on $ {D_{-t}},\ {C_{t}}$ i ${D_{t}}$
s\'on els camins dibuixats en figura, hom ha,
per la mateixa ra\'o  que (\ref{grando}),
$ I_{1t}\rightarrow 0$, $ I_{3t}\rightarrow 0$
per a $ \vert t     \vert\to\infty  $.
Finalment, hem mostrat que (\ref{formal})
\'es una veritable 
soluci\'o de (\ref{quater}), doncs del
problema de Hurwitz.
\end{minipage}
\end{tabular}

\section{Idiosincrasies}
\subsection{La continuaci\'o 
maximal}

Recordem aqu\'\i\ el resultat ben 
conegut que
cada element 
$\left(U,f    \right)$ 
de funci\'o holomorfa
t\'e una continuaci\'o
anal\'\i tica regular
maximal $Q:= \left(S,\pi ,j,F
    \right)  $, on
$S$ \'es una
superf\'\i cie
de Riemann regular 
sobre un conjunt
obert $\Omega$ de $\CI$, amb
una aplicaci\'o de projecci\'o
$p:R\rightarrow\Omega$ que \'es
un biholomorfisme local 
surjectiu,
una 
immersi\'o  holomorfa 
$j\,\colon\, U\rightarrow S$  
tal que $\pi\circ 
j=\id\vert_{U}$ 
i una funci\'o
 holomorfa 
$F\,\colon\, S\rightarrow \ce$
tal que $F\circ j=f$.

Com hem damunt dit, la noci\'o de continuaci\'o anal\'\i tica \'es, en general, incompatible
amb la noci\'o de revestiment topol\`ogic, tot
a que la construcci\'o de l'aplicaci\'o  de projecci\'o de la continaci\'o anal\'\i tica 
en comparteixi 
algunes motivacions inicials.

En aquesta secci\'o
 mostrarem 
la sobreposici\'o 
de punts de frontera 
i de punts interiors de la continuaci\'o, 
anal\'\i tica i, finalment,
la pres\`encia de 
preimatges de discs 
a components
de talla arbitrariament petita.

\subsection{Prolungaments regulars
que no s\'on 
uns
revestiments  }

En aqueste exemple (vegeu \cite{beardon}) mostrem
el comportament de la continuaci\'o 
anal\'\i tica d'un element de funci\'o holomorfa, construit com a un invers local d'un producte de Blaschke.

Sigui $ \DI\subset\CI  $ el disc unitat obert
i $ \{a_n     \}\subset\DI  $ una successi\'o 
de punts de $ \DI  $ que s'acumulen a tot
punts de $\TI:= \partial\DI  $ i
tal que
$\sum_{n=0}^{\infty}(1-\vert a_n\vert)   $
convergeixi.

Amb aquesta hip\`otesi, el producte de Blaschke
\begin{equation}
B(z):=
\prod_{n=1}^{\infty}\frac{z-a_n}{1-\overline{a_n}
z}
\labelle{blas2}
\end{equation}
convergeix uniformement en els conjunts compacts de $ \DI  $ i doncs hi defineix
una funci\'o  anal\'\i tica que s'anula exactament
als punts $ \{a_n     \}  $.

Puix que aquestes zeros s'acumulen a la 
circumfer\`encia unitat $ \TI  $, $ B  $ no pot
ser continuada anal\'\i ticament a algun punt
$ p\in\TI  $: altrament, gracies al teorema
{\bf\ref{successio}} hom hauria $ B\equiv 0  $.

Aix\'\i \ $ \Ti  $ \'es una frontera natural
per $ B  $.

De m\'es, tot factor del producte 
(\ref{blas2}) s\'on automorfismes del disc unitat, doncs 
$$
\vert B(z)\vert =
\prod_{n=1}^{\infty}
\left\vert
\frac{z-a_n}{1-\overline{a_n}
z}
\right\vert
\leq 1,
$$
\'es a dir
$ B(\DI)\subset\DI  $.

Gracies al lema de Schwartz-Pick
(vegeu per exemple \cite{chabat2}, p.348, teorema 2),
$ B  $ \'es
una contracci\'o de la m\`etrica hiperb\`olica de $ \DI  $, doncs
\begin{equation}
0<\Vert B^{\prime}(a_n)  \Vert
\leq
\frac{1-\Vert B^{\prime}(a_n)  \Vert^2}
{1-\vert a_n \vert^2}\leq
\frac{1}{1-\vert a_n \vert^2}
\leq
\frac{1}{1-\vert a_n \vert}.
\labelle{vun2}
\end{equation}

Com, per a tot $ n  $, $ a_n  $ \'es un zero
simple de $ B  $, podem construir una determinaci\'o  $ f_n  $ de $ B^{-1}  $ a $ 0  $
tal que $ f_n (0)=a_n   $.
De m\'es, els $ \{f_n     \}  $ s\'on elements
conectables.
Sigui $ r_n  $ el radi de converg\`encia de 
$ f_n  $ en $ 0  $: llavors $ f_n\left(\DI(0,r_n)    \right) \subset\DI $, si no $ B  $ podria ser
continuada anal\'\i ticament m\'es enll\`a de 
$ \DI  $.

Aix\'\i , $f_n:\DI(0,r_n)\rightarrow 
  f_n\left(\DI(0,r_n)    \right)      $ \'es un
homeomorfisme
i
\begin{equation}
 m\not=n 
\Rightarrow 
 a_m\not\in f_n\left(\DI(0,r_n)    \right).
\labelle{dane}
\end{equation}

Puix que $ f_n  $ \'es conforme en 
$ \DI(0,r_n) $, podem aplicar el teorema
\ref{koebe}, doncs 
$f_n\left(\DI(0,r_n)    \right) \subset
\DI(a_n,r_n
{\vert f_n^{\prime}(0)\vert}/{4}) 
   $.

Gracies a la regla de la cadena, 
$ \vert  f_n^{\prime}(0)B^{\prime}(a_n)   \vert \equiv 1 $.
Construim els $ \{ a_n  \}  $ de manera que
$\vert a_{2n}-
a_{2n+1} \vert\leq (1-a_{2n})^2
   $:
considerant tamb\'e
(\ref{vun2}) i (\ref{dane}),
hom ha, per a tot $n $, 
\begin{equation}
\labelle{tri2}
4\vert a_{2n}-a_{2n+1} \vert \geq r_{2n}
\vert f_{2n}^{\prime}(0)\vert
\geq
r_{2n}(1-\vert a_{2n}\vert),
\end{equation}
aix\'\i 
\begin{equation}
\labelle{quatre2}
r_{2n}\leq
\frac{4\vert a_{2n}-a_{2n+1} \vert }{ (1-\vert a_{2n}\vert)}
\leq 4 (1-\vert a_{2n}\vert),
\end{equation}
\'es a dir
$ r_{2n}\rightarrow 0  $
 per
$ n\rightarrow \infty  $.

Sigui $ \left(S,\pi ,j,F    \right)  $
la 
continuaci\'o anal\'\i tica maximal
d'un $ f_n  $ qualsevol: llavors l'exemple precedent mostra que $ \pi  $ no \'es un 
revestiment topol\`ogic de $ \pi(S)  $, perqu\`e,
sobre qualsevol entorn del punt $ 0\in\DI  $
jeuen discs a radi
arbitrariament petit.

\end{document}